# Engineering design optimisation using Tabu search

**DR ANDY CONNOR** and **DR P JOHN CLARKSON**
Engineering Design Centre, University of Cambridge, Cambridge, UK
**DR SHAHROKH SHAHPAR**
Aerothermal Methods, Rolls-Royce plc, Derby, UK
**DR PAUL LEONARD**
Department of Electrical Engineering, University of Bath, Bath, UK

**ABSTRACT**

This paper describes an optimisation methodology that has been specifically developed for engineering design problems. The methodology is based on a Tabu search (TS) algorithm that has been shown to find high quality solutions with a relatively low number of objective function evaluations. Whilst the methodology was originally intended for a small range of design problems it has since been successfully applied to problems from different domains with no alteration to the underlying method. This paper describes the method and it's application to three different problems. The first is from the field of structural design, the second relates to the design of electromagnetic pole shapes and the third involves the design of turbomachinery blades.

## 1 INTRODUCTION

Numerical optimisation techniques can be applied as a tool in the design of engineering systems and are particularly applicable to parametric design problems, where the general form or type of solution is known but it is necessary to determine values for the design variables to ensure that the system produces the desired response.

In this paper three different examples of parametric design are considered. In the first example a ten bar truss is designed in order to achieve a particular structural performance. In the second example, the shape of an electromagnetic pole is determined such that a given electrical field is produced. In the final example, a turbine blade is designed using computational models developed at Rolls Royce. This paper is not intended to be either a complete introduction to TS or a fully verified application paper. The intention is to illustrate that the TS implementation can be successfully applied to a range of problems without specific modification and produces high quality results in all cases.

## 2   TABU SEARCH

TS (1) is a metaheuristic which is used to guide optimisation algorithms in the search for a globally optimal solution. The algorithm uses flexible memory cycles of differing time spans to force the search out of local optima and to provide strategic control to progress the search through the solution space. All the results in this paper have been obtained using an implementation that was originally designed for the optimisation of fluid power circuits (2).

### 2.1   Short Term Memory
The most simple implementation of TS is based around the use of a hill climbing algorithm. Once the method has located a locally optimal solution the use of the short term memory, or tabu restrictions, ensures that the search does not return to the optimum after the algorithm forces the search out in a new direction. In the TS implementation used in this work, the short term memory contains a list of the last $n$ visited solutions and these are classed as tabu.

The effect of this concept can be illustrated by considering the diagram shown in Figure 1. This shows a contour plot of a two-dimensional function which contains one local and one global optimum and the aim of the search is to find the location with the lowest value.

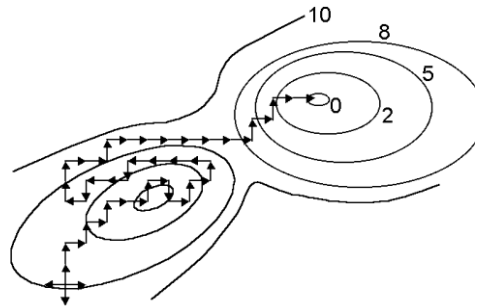

**Fig. 1  Tabu search vector**

From the indicated start position the local search algorithm quickly locates the locally optimal. When the search reaches the local optimum the aggressive nature of the TS forces the algorithm out of the optimum in the direction that increases the objective function by the smallest amount. Because the last $n$ visited solutions are classed as tabu, the search cannot leave the optimum along the reverse trajectory from which it entered and once it has left the optimum it cannot enter it again. The algorithm therefore forces the search to climb out of the local optimum and in due course it successfully locates the global optimum.

### 2.2   Search Intensification and Diversification
The TS short term memory enables the method to leave locally optimal solutions in the quest for the global optimum of a function. However, short term memory alone does not ensure that the search will be both efficient and effective. Search intensification and diversification techniques are often used first to focus the search in particular areas and then to expand the search to new areas of the solution space. This is normally achieved by the use of longer term memory cycles.

Intermediate and long term memory cycles generally use similar lists of previously visited solutions to guide the search. In the specific implementation used in this work the intermediate term memory cycle is based on a list of the *m* best solutions found so far. This list is therefore only updated when a new improved solution is found as opposed to whenever a move is made. At certain stages throughout the search process a degree of intensification is achieved by reinitialising the search at a new point generated by considering similarities between the solutions contained in the intermediate memory list. In the implementation used in this paper, diversification is achieved by using a simple random refreshment although more strategic diversification could be implemented through the use of long term memory.

## 3    STRUCTURAL DESIGN

The structural problem considered in this paper is the optimisation of the ten bar truss. In the standard approach to this problem the spatial layout of the truss is constrained as shown in Figure 2 by fixing the position of the nodes and TS has been applied a number of times to this problem (3,4). The truss is assumed to consist of an idealised set of pin jointed bars connected together at the nodes. The design optimisation problem is to find the cross sectional areas of each member such that the mass is minimised.

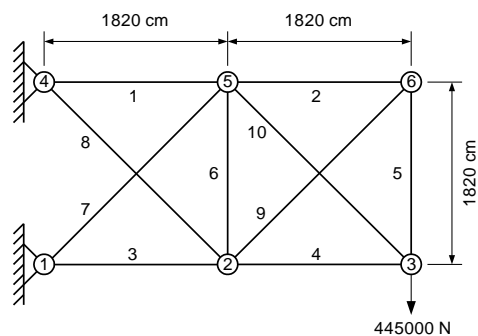

**Fig. 2 Ten bar truss**

The problem has been expanded in this paper to allow the spatial layout of the truss to be adjusted as well as the cross sectional area of the members. The built in nodes are fixed in position, as is the loaded node. This variation of the problem therefore introduces six new parameters that determine the position of nodes 2, 5 and 6.

Essentially, all of the design parameters are continuous. However, the nodal positions have a minimum allowable change of 1cm whilst for the cross sectional areas the minimum allowable change is $0.01 cm^2$. In reality, these cross sectional areas would be limited to discrete values corresponding to available stock material. The material for the truss is aluminium with Young's modulus of $6.88 \times 10^6$ $N/cm^2$ and material density is $2.7 \times 10^{-3}$ $kg/cm^3$. Each member is modelled as a solid circular cross section.

Despite the simplicity of the ten bar truss example, it is still a reasonably constrained problem due to the difficulty in finding high quality, i.e. low mass, solutions that do not violate either the buckling or stress constraints. The constraints on the problem are that each member should not have a stress that exceeds 17,200 N/cm$^2$, buckle under Euler buckling criteria, and have a length less than 15cm. Table 1 shows the numeric values for the x,y coordinates of the node positions of a typical solution expressed relative to the lower fixed node. Table 2 shows the cross sectional areas of the members where some can be seen to be minimum area members. This solution has a mass of 1598kg and was found in 12004 evaluations.

Table 1. Node coordinates

| Node | Coordinates cm) |
|---|---|
| $x_2,y_2$ | 445,-61 |
| $x_5,y_5$ | 807,408 |
| $x_6,y_6$ | 1197,-112 |

Table 2. Member areas

| Member | Area (cm$^2$) |
|---|---|
| $A_1$ | 60.39 |
| $A_2$ | 16.6 |
| $A_3$ | 183.17 |
| $A_4$ | 0.01 |
| $A_5$ | 239.9 |
| $A_6$ | 3.04 |
| $A_7$ | 0.01 |
| $A_8$ | 1.42 |
| $A_9$ | 310.26 |
| $A_{10}$ | 47.9 |

Such minimum area members can be removed from the structure to produce a reduced topology solution provided that the change does not produce a violation of constraints or so great a reduction in topology that the structure begins to act as a mechanism. Simple violations may be adjusted out by the designer, however it is important to realise that if a minimum area member is transferring significant force to a node, then removing that node may result in significant changes in the response of other members which become increasingly difficult to gauge. The solution in Figure 3 resulted from the tabu search method. Removing the minimum area members produces the solution shown in Figure 4. This reduced topology has no constraint violations and thus needs no further adjustment.

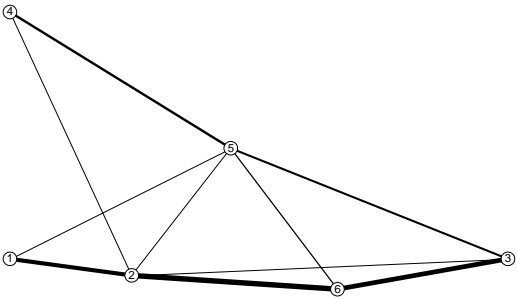
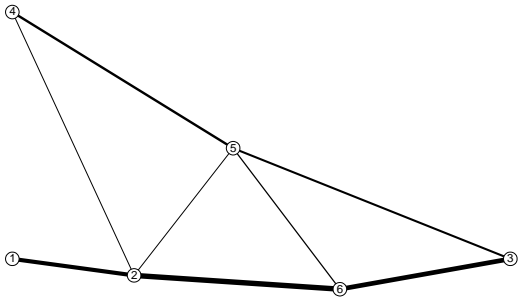

Fig. 3 Truss solution                Fig. 4 Adjusted truss

These solutions both exhibit considerably lower mass than can be achieved by using a standard ten bar truss. The solution also has smooth load transmission paths and a more elegant appearance. Whilst the number of evaluations is relatively high, comparisons with other methods show that the required number of evaluations is less than for methods such as simulated annealing.

## 4 MAGNETIC POLE SHAPE DESIGN

MRI imaging applications require a high field homogeneity within a given region. This uniformity of field is often achieved by utilising iron poles with a particular shape (5). The problem in this paper involves a pole shape problem. A pair of rotationally symmetric poles are used to shape the field as shown in Figure 5. For simplicity the source of the field is not modelled and the poles are assumed to be driven by auxiliary coils or permanent magnets.

A number of different pole parameterisations have been investigated. Figure 6 illustrates a pole described by four parameters that has two ramped sections. The other parameterisations were specified so that solutions had one, three and four ramps. A parameterisation that had no structure but used the heights of sixteen points along the length of the pole is not discussed in this paper.

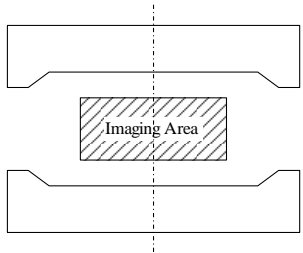
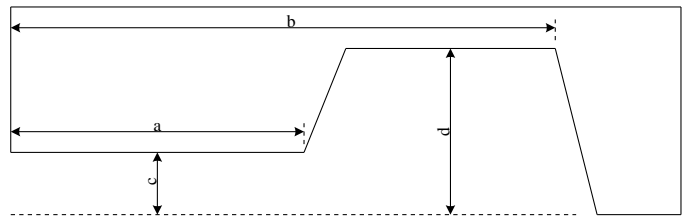

**Fig. 5 MRI application**  **Fig. 6 Four parameter pole shape**

Typical pole shapes for each of the parameterisations are shown in Figure 7, where the sub-plots (a) to (d) show solutions from one ramp to four ramps consecutively.

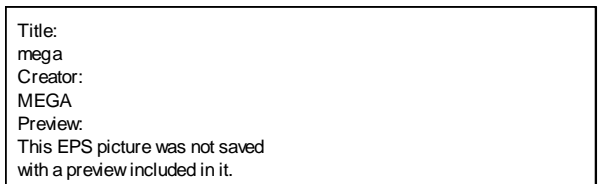
(a)
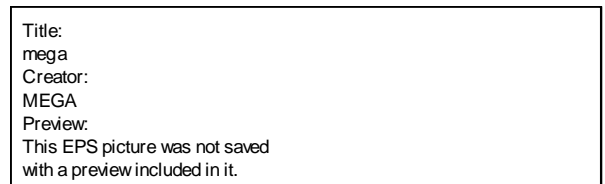
(b)
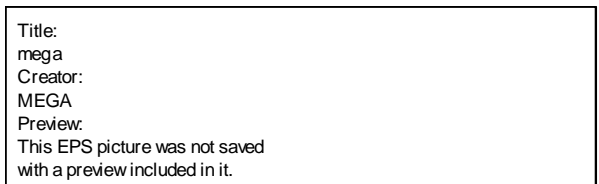
(c)
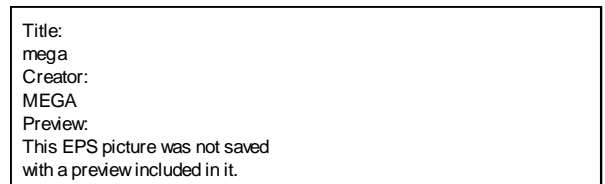
(d)

**Fig. 7 Typical solutions**

Table 3 shows the range of objective function values, which is an indication of the uniformity of the field, for a number of runs. The values in the table show the maximum and minimum values found after 20, 100 and 1000 evaluations as well as for the converged search.

**Table 3 Results**

|        | Objective Function |     |      |           |
|--------|--------------------|-----|------|-----------|
| Scheme | 20                 | 100 | 1000 | Converged |

| 1 | 1006-3897 | 597-1118 | 520-595 | 520-595 |
| 2 | 610-2252 | 284-1136 | 126-250 | 112-149 |
| 3 | 1105-4474 | 213-2343 | 66-425 | 55-217 |
| 4 | 6065-13256 | 576-1419 | 183-429 | 98-422 |

Due to the probabilistic nature of the algorithm, including random start points, the spread of solution quality is quite high. However, it can be seen that the method is rapidly improving the initial solutions and in the case of the two ramp scheme is finding solutions within a fairly narrow band of objective function values. The three ramp scheme has located a solution with a lower objective function value but the standard deviation of objective function values is considerably higher.

## 5   TURBINE BLADE DESIGN

In this section, the application of the Tabu search in the design of the turbomachinery blades is discussed. The objective of the design effort is to minimize the secondary-flow kinetic energy (SKE) of a nozzle-guide vane (NGV). Due to high loading i.e. large flow turning and low aspect ratio, the secondary flows are the dominant features of the NGV studied here. The reduction in the SKE of the blade leads to a reduction in loss and hence higher efficiencies and ultimately improved specific fuel consumption for the whole engine can be obtained.

Since the efficiencies of the blades are already very high, additional improvements can only come from three-dimensional design of blades. The three dimensional aerodynamic design tool used here is based on the FAITH (Forward And Inverse THree-dimensional) linear design system. The FAITH system had successfully been applied in the design of Turbines (6) and Compressor Blades (7). The forward design mode of FAITH has been used as a fast CFD code requiring no more than a few minutes of computation time.

A three-dimensional pressure correction based solver is used to perform the CFD calculations for base and the perturbations. The method has a low level of numerical viscosity due to its control volume upwinding algorithm. The code has been shown to give very accurate predictions of the static pressure field for turbomachinery applications and good qualitative predictions of the viscous phenomena and secondary flows on relatively coarse meshes.

In this paper, results for two design categories namely lean (YCEN) and sweep (XCEN) are presented. The parameters are specified at seven control points along the span. All the intervening mesh points are linearly interpolated. Typical perturbations for YCEN and XCEN design categories at mid height are shown in figures 8 and 9 respectively.

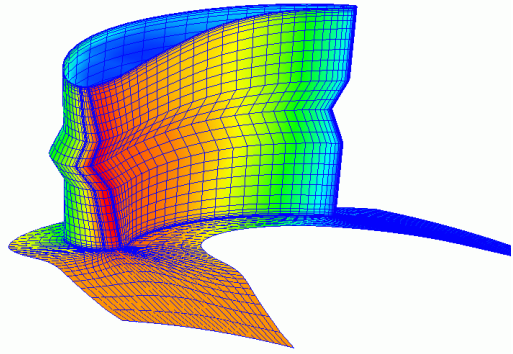 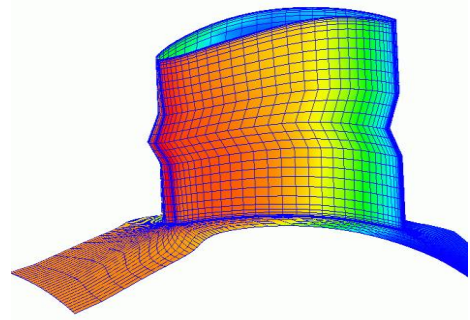

**Fig. 8  Circumferential view of the YCEN perturbation at mid-height**     **Fig. 9  Meridional view of the XCEN perturbation at mid-height**

A summary of the results is given in table 4. The best reduction in cost function is achieved by combining YCEN and XCEN together. Although, there has been a rise in SKE of the order of  ~20% due to non-linearity effects, still a ~40% reduction in SKE according to non-linear results have been achieved.  The small numbers quoted here for the SKE is due to the non-dimensionalisation of the secondary-kinetic energy by the mass-mean dynamic head at the exit-plane on the NGV.

**Table 4  Comparison of reductions in the cost function**

|              | SKE (LIN) | SKE (NLIN) |
|--------------|-----------|------------|
| BASE         | -------   | 0.04466    |
| XCEN         | 0.02970   | 0.03964    |
| YCEN         | 0.02570   | 0.03058    |
| XCEN + YCEN  | 0.01767   | 0.02704    |

Figures 10 and 11 show the contours of the SKE at the exit plane of the base and optimised geometry using 14 design parameters including lean and sweep design categories respectively.

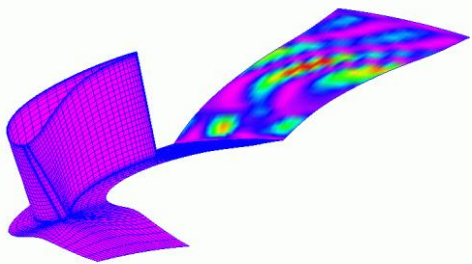 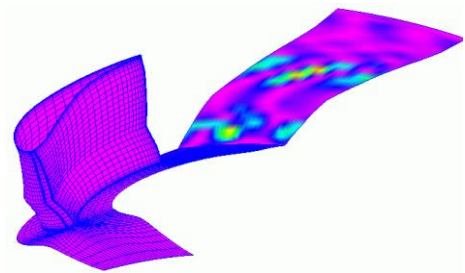

**Fig. 10  Contours of the secondary-kinetic energy at the exit plane of the base geometry (non linear results)**     **Fig. 11  Contours of the secondary-kinetic energy at the exit plane of the optimised geometry (non linear results)**

Figure 12 show the mass-mean plots of SKE at the exit-plane of the vane and the reductions achieved by the optimiser using the linear BFAITH code and verified by the non-linear calculations.

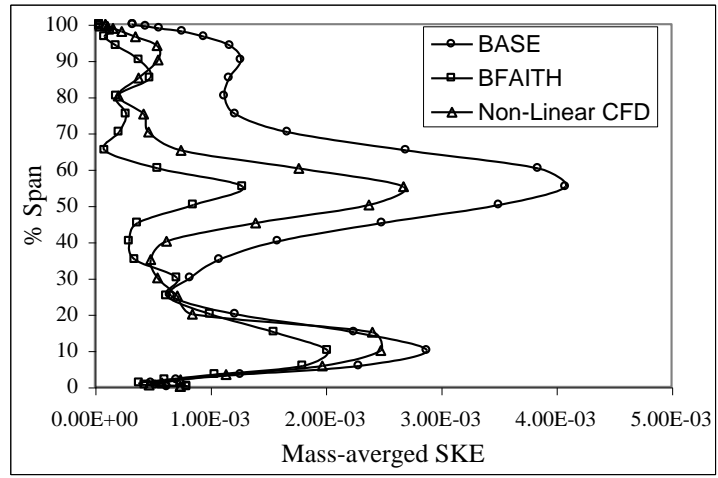

**Fig. 12 Spanwise variation of the circumferentially mass-averaged SKE at the exit-plane of the base, and optimised geometry.**

## 6 CONCLUSIONS

This paper has presented a brief introduction to the Tabu search method and has provided a small sample of application problems to illustrate that the implementation developed can be applied to a wide range of problems. The sample problems represent different domains, where each problem has it's own representation that has been used to define the optimisation. Given that detailed results have not been presented it is difficult to attempt to claim any conclusions other than the generality of the search algorithm. However, future work and publications will show that in all cases the method has found good design solutions with no modification of the underlying algorithm.

**ACKNOWLEDGEMENTS**

The authors would like to thank Rolls Royce plc for permission to publish the results regarding the application of the optimisation method to turbomachinery blade design.